\newtheorem{theorem}{Theorem}
\newtheorem{corollary}{Corollary}

\newtheorem{lemma}{Lemma}
\newtheorem{remark}{Remark}
\newtheorem{example}{Example}

\documentclass[12pt]{article}
\usepackage{amsfonts}

\usepackage{amsmath}
\usepackage{graphicx}

\newcommand{\R}{\mathbb{R}}
\newcommand{\Q}{\mathbb{Q}}

\newcommand{\veps}{\varepsilon}
\newcommand{\ve}{\varepsilon}

\newcommand{\halmos}{$\sqcup\!\!\!\!\sqcap$}

\newcommand{\rmd}{{\rm d}}
\newcommand{\rmi}{{\rm i}}

\newcommand{\dto}{\downarrow}

\newcommand{\be}{\begin{equation}}
\newcommand{\ee}{\end{equation}}
\newcommand{\ben}{\begin{equation*}}
\newcommand{\een}{\end{equation*}}
\newcommand{\ba}{\begin{aligned}}
\newcommand{\ea}{\end{aligned}}

\newcommand{\ga}{\alpha}
\newcommand{\gb}{\beta}
\newcommand{\gk}{\kappa}

\newcommand{\wh}{\widehat}

\newcommand{\whh}{\widehat{h}}
\newcommand{\whH}{\widehat{H}}
\newcommand{\whV}{\widehat{V}}
\newcommand{\whL}{\widehat{L}}

\newcommand{\wht}{\widehat{t}}
\newcommand{\whU}{\widehat{U}}
\newcommand{\whX}{\widehat{X}}

\newcommand{\Xbar}{\overline{X}}

\newcommand{\Gtau}{G_{\tau_u-}}

 \newcommand{\tk}{\kappa_{Z,Y}}

\numberwithin{equation}{section}
\numberwithin{theorem}{section}
\numberwithin{lemma}{section}
\numberwithin{proposition}{section}
\numberwithin{remark}{section}
\numberwithin{corollary}{section}
\numberwithin{example}{section}
\begin{document}

\title{The time at which a  L\'evy
process creeps
}
\author{Philip S. Griffin and Ross A. Maller\thanks{Research partially supported by
ARC Grant DP1092502}
\\Syracuse University and Australian National University}
\maketitle

\begin{abstract}
We show that if a L\'evy
process creeps then,  as a function of $u$, the renewal function $V(t,u)$ of the bivariate ascending ladder process $(L^{-1},H)$ is absolutely continuous on $[0,\infty)$ and left differentiable on $(0,\infty)$, and the left derivative at $u$ is proportional to the (improper) distribution function of the time at which the process creeps over level $u$, where the constant of proportionality is $\rmd_H^{-1}$, the reciprocal of the (positive) drift of $H$.
This yields the (missing) term due to creeping in the recent quintuple law of Doney and Kyprianou (2006).
As an application, we
derive a Laplace transform identity
which generalises the second factorization identity.  We also relate Doney and Kyprianou's extension of Vigon's \'equation amicale invers\'ee to creeping.
Some results concerning the ladder process of $X$, including the
second factorization identity, continue
to hold for a general bivariate subordinator, and are given in this
generality.
\end{abstract}




\vspace*{10pt} \setcounter{equation}{0}
\section{Introduction}\label{s1}
Let $X=\{X_{t}: t \geq 0 \}$,
$X_0=0$,  be a real-valued  L\'{e}vy process
with characteristic triplet $(\gamma, \sigma^2, \Pi_X)$,
thus the characteristic function of $X$
is given by the L\'{e}vy-Khintchine
representation, $Ee^{i\theta X_{t}} = e^{t \Psi_X(\theta)}$,
where
\begin{equation}\label{lrep}
\Psi_X(\theta) =
 \rmi\theta \gamma - \sigma^2\theta^2/2+
\int_{\R}(e^{\rmi\theta x}-1-
\rmi\theta x \mathbf{1}_{\{|x|<1\}})\Pi_X(\rmd x),
\ {\rm for}\  \theta \in \mathbb{R}, \ t\ge 0.
\end{equation}
$X$ is said to {\it creep across a level} $u>0$ if
$P(\tau_u<\infty, X_{\tau_u}=u)>0$
where
\ben
\tau_u= \inf \{t \geq0: X_t>u \}.
\een
Our initial interest is in the {\it time at which  $X$ creeps}.  Thus we introduce  the (improper) distribution function
\begin{equation}\label{ct2}
p(t,u)=P(\tau_u\le t, X_{\tau_u}=u),\ u>0, t\ge 0.
\end{equation}
We prove certain regularity properties of  $p(t,u)$ which allow us to relate it to the renewal function of the bivariate ascending ladder process of $X$.  This yields the missing term, due to creeping, in Doney and
Kyprianou's (2006) {\it quintuple law}.
Using the quintuple law, we derive a Laplace transform identity which generalises
the {\em second factorization identity}
due to Percheskii and Rogozin (1969).   We also relate creeping to  Doney and Kyprianou's extension of the \'equation amicale invers\'ee of Vigon (2002).
Some of these results  extend from the bivariate ladder process to general bivariate subordinators, and we develop several of the results in this setting.  In particular, it appears to have gone previously  unnoticed, that the   second factorization identity is a special case of a general transform result for bivariate subordinators.
The results in the fluctuation setting are stated in Section \ref{s3}, with their proofs given in Sections \ref{s5}.  The general bivariate subordinator case is developed in Section \ref{s4}.

By a  compound Poisson process we will mean a
L\'evy process with finite L\'evy measure, no Brownian component and zero drift.
The indicator of an event $A$ will be denoted by
${\bf 1}_A$, or sometimes by ${\bf 1}(A)$,
and we adopt the convention that the inf of the empty set is $+\infty$.


\setcounter{equation}{0} \section{Fluctuation Setup}\label{s2}

We need some notation, which is very standard in the area.
Let $(L_s)_{s \geq 0}$ be the local time at the maximum, and
$(L^{-1}_s,H_s)_{s \geq 0}$
the weakly\footnote{The distinction between weak and
strict only makes a difference when $X$ is compound Poisson.}
ascending bivariate ladder process of $X$.
Bertoin (1996), Chapter VI, and Kyprianou (2006), Chapter 6,
give detailed discussions of these processes
and their properties; see also
Doney (2005).
When $X_t\to -\infty$ a.s., $(L^{-1},H)$ is defective and may
be obtained from  a nondefective bivariate subordinator
$({\cal L}^{-1},{\cal H})$ by exponential killing at some
rate $q > 0$, say.
Thus
\be\label{H}
(L^{-1}_s, H_s)=\begin{cases} ({\cal L}^{-1}_s, {\cal H}_s) & \text{ if } s<e(q),\\
(\infty,\infty) & \text{ if } s\ge e(q),
\end{cases}\ee
where $e(q)$ is independent of $({\cal L}^{-1},{\cal H})$ and has exponential
distribution with parameter $q$.
In the case that
$(L^{-1},H)$ is nondefective there is no need to introduce
exponential killing and we set
$(L^{-1},H)=({\cal L}^{-1},{\cal H})$.

We denote the bivariate L\'{e}vy measure of
$({\cal L}^{-1},{\cal H})$
by $\Pi_{{L}^{-1},{H}}(\cdot,\cdot)$, and its marginals by  $\Pi_{{L}^{-1}}(\cdot)$ and
$\Pi_{{ H}}(\cdot)$.
The Laplace exponent $\kappa(a,b)$
of  $(L^{-1},H)$ is given by
\begin{equation} \label{kapdef}
e^{-\kappa(a,b)}  = E (e^{-a{L}^{-1}_1 -b{H}_1};1<L_\infty) =e^{-q} E e^{-a{\cal L}^{-1}_1 -b{\cal H}_1}
 \end{equation}
for values of $a,b\in \R$ for which the expectation is finite.  It
may be written
\begin{eqnarray} \label{kapexp}
\kappa(a,b)
&=&
q+\rmd_{L^{-1}}a+\rmd_Hb+\int_{t\ge0}\int_{x\ge0}
\left(1-e^{-at-bx}\right)
\Pi_{{L}^{-1}, {H}}(\rmd t, \rmd x),
\end{eqnarray}
where $\rmd_{L^{-1}}\ge 0$ and $\rmd_H\ge 0$ are drift constants.
The bivariate renewal function of
$(L^{-1},H)$ is\footnote{Throughout, we will
write $\Pi_{{L}^{-1}, {H}}(\cdot,\cdot)$ and
$V(\cdot,\cdot)$ with the
time variable in first position, followed by the space variable.
This is at variance with some established literature,
but seems desirable for consistency.}
\begin{eqnarray}\label{Vkdef}
V(t,x)=\int_0^\infty P(L_s^{-1}\le t,H_s\le x)\rmd s
=\int_0^\infty e^{-qs}P({\cal L}^{-1}_s \le t,{\cal H}_s \le x)\rmd s.
\end{eqnarray}
It has Laplace transform
\begin{equation}\ba\label{Vkap}
\int_{t\ge 0}\int_{x\ge 0}e^{-at-bx} V(\rmd t,\rmd x)
=\frac{1}{\kappa(a,b)}
\ea\end{equation}
for all $a, b$ such that $\kappa(a,b)>0$. The positivity condition on $\gk$ clearly holds when $a,b\ge 0$ and either $a\vee b>0$ or $q>0$.

Let $\whX_t=-X_t$, $t\ge 0$ denote the dual process, and $(\whL^{-1}, \whH)$ the corresponding {\it strictly} ascending bivariate
ladder processes of $\whX$.  This is the same as the weakly ascending process if $\whX$ is not compound Poisson. The
definition of $(\whL^{-1}, \whH)$ when $\whX$ is compound Poisson is as the limit of the  ascending bivariate ladder
process of $\whX_t-\ve t$ as $\ve\downarrow 0$.
All quantities relating to $\whX$ will be denoted in the obvious way; for example
$\Pi_{{\whL}^{-1}, {\whH}}(\cdot,\cdot)$, $\widehat \kappa(\cdot,\cdot)$ and $\whV(\cdot,\cdot)$.
We choose the normalisation of the local times $L$ and $\whL$ so that
the Weiner-Hopf factorisation  takes the form
\begin{equation}\label{WH}
\kappa(a,0) \wh\kappa(a,0)
= a,\ \ a\ge0.
\end{equation}
This would not be possible  in the compound Poisson case if $({\whL}^{-1}, {\whH})$ were the weak bivariate ladder process; see Section 6.4 of Kyprianou (2006).


\setcounter{equation}{0} \section{Creeping Time}\label{s3}

It is well known that  $X$ creeps across some $u>0$ iff
$X$ creeps across all $u>0$, in which case we say that $X$ creeps.  A necessary and sufficient condition for creeping is that $\rmd_H>0$; see Theorem VI.19 of Bertoin (1996).
Our first result describes the (improper) distribution function of the time at which $X$ creeps across level $u$.

\begin{theorem}[Creeping Time] \label{thmct1}\
\newline \noindent (i)\
The following are equivalent:
\be\label{C1}
\text{ $p(t,u)>0$ for some $t>0$, $u>0$; }
\ee
\be \label{C2}
\text{  $p(t,u)>0$ for all $u>0$ and all
$t$ sufficiently large (depending on $u$);}
\ee
\be\label{C3}
\text{ $p(t,u)>0$ for all $t>0$ and all
$u$ sufficiently small (depending on $t$);}
\ee
\be\label{C4}
 \text{ $\rmd_H>0$.}
\ee
\noindent
(ii)\
If $\rmd_H>0$, then for every $t\ge 0$, $V(t,0)=0$,
$V(t,\cdot)$ is absolutely continuous on $[0,\infty)$ with a left continuous left derivative on $(0,\infty)$, and  for each $u\in (0,\infty)$ satisfies
\begin{equation}\label{ct1}
p(t,u)
=\rmd_H \frac{\partial_-}{\partial_- u}V(t,u)
\end{equation}
where ${\partial_-}/{\partial_- u}$ denotes the left partial derivative in $u$.
\newline\noindent
(iii)\
If $\rmd_H>0$ and  $X$ is not compound Poisson with positive drift,
then $V(t,\cdot)$ is differentiable and  $p(t,\cdot)$ is
continuous on  $(0,\infty)$ for each $t\ge 0$,
and $p(\cdot,u)$ is  continuous on $[0,\infty)$ for each $u>0$.
\end{theorem}

\begin{remark}\label{R1}
{\rm
(i)\
Unlike the creeping case, it is possible that $X$ creeps over some $u$ but not all $u$ by a fixed time $t>0$.  This is illustrated in Examples
\ref{E2} and \ref{E3}.  Nevertheless, from Theorem  \ref{thmct1},
$X$ creeps over some $u>0$ by some time $t>0$  iff $X$ creeps,
and we obtain the generalisation \eqref{ct1} of
Kesten and Neveu's formula for the probability of eventually creeping over $u$; see pp.\ 119--121 of Kesten (1969).

(ii)\
It follows immediately from Theorem \ref{thmct1} that
\be\label{Vdens}
\rmd_H V(\rmd t, \rmd u)= P(\tau_u\in \rmd t, X_{\tau_u}=u)\rmd u.
\ee
This formula has already been noted by
Savov and Winkel (2010), p.8, and attributed to Andreas Kyprianou.
Conversely, from Savov and Winkel's observation \eqref{Vdens}, it follows that $\rmd_H V(t, \rmd u)$ has a density given by $p(t,u)$ for a.e.$\ u$.
This however gives no information about $p(t,u)$ for a given level $u$.  One of the main points of Theorem
\ref{thmct1} is that $p(t,u)$ \emph {is
the left derivative of $\rmd_H V(t, u)$ for every}  $u>0$, $t\ge 0$.  This is particularly relevant in the quintuple law below.

(iii)
In the case that $X$ is a subordinator which creeps, the Laplace transform of the time at which it creeps over $u$ is given by
\be\label{rd}
E(e^{-\ga \tau_u};X_{\tau_u}=u)=\rmd_X v^\ga(u)
\ee
where $v^\ga$ is the bounded continuous density of the resolvent kernel
\ben
V^\ga(\rmd u)=\int_0^\infty e^{-\ga t} P(X_t\in \rmd u)\rmd t;
\een
see page 80 of Bertoin (1996). This in principle gives the distribution of the time at which $X$ creeps over $u$.  Indeed
\ben\ba
v^\ga(u)&=\frac{\rmd}{\rmd u}\int_0^\infty e^{-\ga t} P(X_t\le u)\rmd t
&= \frac{\rmd}{\rmd u}\int_0^\infty\ga e^{-\ga t}\rmd t \int_0^tP(X_s\le u)\rmd s.
\ea\een
Hence if the derivative could be moved inside the integral, from \eqref{rd} we would obtain
\be\label{md}
P(\tau_u\le t, X_{\tau_u}=u)= \rmd_X\frac{\partial}{\partial u} \int_0^tP(X_s\le u)\rmd s.
\ee
Since
\ben
\rmd_HV(t,u)=\rmd_X\int_0^t P(X_s\le u)\rmd s
\een
when $X$ is a subordinator, it follows from \eqref{ct1} that
\ben
P(\tau_u\le t, X_{\tau_u}=u)= \rmd_X\frac{\partial_-}{\partial_- u} \int_0^tP(X_s\le u)\rmd s.
\een
Thus \eqref{md} is correct  provided ${\partial}/{\partial u}$ is replaced by ${\partial_-}/{\partial_- u}$.  Conversely one can use \eqref{ct1} to give an alternative proof of \eqref{rd}.

(iv)\
Theorem \ref{thmct1} is concerned with regularity of $V(t,\cdot)$.  Some information about regularity of $V(\cdot,u)$ may be gleaned from Theorem 5 of Alili and Chaumont (2001), from which
it follows that
$V(\cdot,u)$ is absolutely continuous for each $u>0$ provided 0 is
regular for both $(-\infty, 0)$ and $(0, \infty)$, for $X$.

}
\end{remark}
\bigskip

The quintuple law is a fluctuation identity, due to Doney and
Kyprianou (2006), describing the joint distribution of five random variables associated with the first passage of $X$ over a fixed level $u>0$ when $X_{\tau_u}>u$.  Using Theorem \ref{thmct1},
we are  able to account for the contribution due to creeping, that is the term when $X_{\tau_u}=u$. To give the result
set
\[
\Xbar_{t} = \sup_{0\le s\leq t}X_s\quad {\rm and }\quad G_t= \sup\{0\le s\le t: X_s= \Xbar_{s}\}.
\]
The quintuple law concerns the following  quantities:
\vskip .1in
\smallskip$\bullet$ First Passage Time Above Level $u$:
$ \tau_u= \inf \{t \geq0: X_t>u \}$;

\smallskip$\bullet$ Time of Last Maximum Before Passage: $G_{\tau_u-}$;


\smallskip$\bullet$ Overshoot Above Level $u$: $X_{\tau_u}-u$;

\smallskip$\bullet$ Undershoot of  Level $u$: $u-X_{\tau_u-}$;

\smallskip$\bullet$ Undershoot of the Last Maximum Before Passage:
$u-\overline X_{\tau_u-}$.
\vskip .1in

\begin{theorem}[Quintuple Law with Creeping]\label{thmct2}\
Fix $u>0$;
then for $ x\ge 0$, $v\ge 0$, $0\le y\leq u \wedge v$, $s \geq0$
and $t\geq0$
\be\ba\label{DKc}
&P\left(X_{\tau_u}-u \in \rmd x,
u-X_{\tau_u-} \in \rmd v,
u-\overline{X}_{\tau_u-} \in \rmd y,
\tau_u-\Gtau \in \rmd s, \Gtau \in \rmd t\right)
\\
&\ \
={\bf 1}_{\{x>0\}}
|V(\rmd t,u-\rmd y)|\widehat{V}(\rmd s, \rmd v-y) \Pi_{X}(\rmd x+v)
+\rmd_H \frac{\partial_-}{\partial_- u}V(\rmd t,u)
\delta_0(\rmd s,\rmd  x,\rmd v,\rmd y),
\ea\ee
where $\delta_0$ is a point  mass at the origin,
and with the convention that the term containing the differential
${\partial_-}V(\rmd t,u)/{\partial_- u}$ is absent when
$ \rmd_H=0$ (in which case  ${\partial_-}V(t,u)/{\partial_- u}$
need not be defined).
\end{theorem}

The contribution to \eqref{DKc} for $x>0$ is Doney and Kyprianou's quintuple law.  Theorem \ref{thmct2} then
follows easily for a.e. $u$ from Savov and Winkel's observation \eqref{Vdens}, but this is clearly unsatisfactory, since it says
nothing about any given $u$.  To get the result for every $u$, \eqref{ct1} is needed.
As a simple consequence of Theorem \ref{thmct2}, we record the joint distribution of the first passage time and overshoot of a level $u>0$;
\begin{corollary}\label{jtop}  Fix $u>0$.  Then for $x,r\ge 0$
\ben\ba
P(X_{\tau_u}-u \in \rmd x, \tau_u\in \rmd r)
&=I(x>0)\int_{0\le s\le r}\int_{0\le y\le u}|V(\rmd s,u-\rmd y)|\Pi_{{ L}^{-1}, { H}}(\rmd r-s,y+\rmd x)\\
&\qquad +\rmd_{ H}\frac{\partial_-}{\partial_- u}V(\rmd r,u)\delta_{\{0\}}(\rmd x).
\ea\een
\end{corollary}
In particular the distribution of the first passage time is
\ben\label{pto2}
P(\tau_u\in \rmd r)
=\int_{0\le s\le r}\int_{0\le y\le u}|V(\rmd s,u-\rmd y)|\Pi_{{ L}^{-1}, { H}}(\rmd r-s,(y,\infty))
+\rmd_{ H}\frac{\partial_-}{\partial_- u}V(\rmd r,u).
\een

Using the quintuple law, Doney and Kyprianou (2006) (Corollary 6) obtain the following useful extension of
the \'equation amicale invers\'ee of Vigon (2002);
for $s\ge 0$ and $x>0$,
\be\label{DKcor6}
\Pi_{L^{-1},H}(\rmd s,\rmd x)
=\int_{v\ge 0} \whV(\rmd s, \rmd v)\Pi_X(\rmd x+v).
\ee
They
state this result for $s>0, x>0$, but their proof works equally well
when $s=0$.
Observe that \eqref{DKcor6} gives no information about
$\Pi_{L^{-1},H}$ on $\{(s,0):s>0\}$.
When $X$ is not compound Poisson, $\Pi_{L^{-1},H}(\rmd s, \{0\})$ is
seen to relate to creeping as the following result shows:

\begin{theorem}\label{creepPi}
Assume $X$ is not compound Poisson. Then
$X$ creeps iff $\Pi_{L^{-1},H}(\rmd s, \{0\})$ is not the zero measure.
\end{theorem}

Despite the connection with creeping, it is easily seen that the jumps
of $(L^{-1}, H)$ for which $\Delta H=0$ do not occur when $X$ creeps
over a fixed level; see \eqref{ncreepPi1} in Section \ref{s5}.
As a consequence of Theorem \ref{creepPi} we are  able to
characterise when \eqref{DKcor6} holds for all $s\ge 0$ and $x\ge 0$:

\begin{theorem}\label{CDK}
\eqref{DKcor6} holds for all  $s\ge 0$, $x\ge 0$ iff $X$ does not creep.
\end{theorem}
\bigskip

When $X$ is compound Poisson it does not creep,
and we can deduce from Theorem \ref{CDK} that
\be\label{CPV2}
\Pi_{L^{-1},H}(\rmd s, \{0\})=
\int_{v\ge 0} \whV(\rmd s, \rmd v)\Pi_X(\{v\}), \ s\ge 0.
\ee
If $\Pi_X$ is diffuse then $ \Pi_{L^{-1},H}(\rmd s, \{0\})$
reduces to the zero measure, but in general it may have positive mass.
Thus  Theorem \ref{creepPi} cannot be extended to the compound
Poisson case.

The next result is an application of  Theorem \ref{thmct2}
to computing a quadruple Laplace transform.  The finiteness conditions on $\gk$, below, clearly hold when $\mu, \rho, \lambda, \nu, \theta\ge0$, and in that case $\gk(\nu,\mu)> 0$ except when $\nu=\mu=0$ and  $q=0$ in \eqref{kapexp}. More generally the conditions  allow for distributions with exponential moments, which can arise quite frequently in applications.

\begin{theorem}[A Laplace Transform Identity]\label{thmct3}\
\newline
Fix $\mu, \rho, \lambda, \nu, \theta$ so that $\gk(\theta,\mu+\lambda), \gk(\theta,\rho)$ are finite and $\gk(\nu,\mu)>0$.

\noindent
(i)\
If $\lambda\ne \rho-\mu$ then
\begin{equation}\label{e3}
\int_{u\ge 0}e^{-\mu u}
E\left(e^{-\rho(X_{\tau_u}-u)-\lambda(u-\Xbar_{\tau_u-})
-\nu \Gtau-\theta(\tau_u-\Gtau)}; \tau_u<\infty\right)\rmd u
=\frac {\kappa(\theta,\mu+\lambda)-\kappa(\theta,\rho)}
{(\mu+\lambda-\rho)\kappa(\nu,\mu)}.
\end{equation}
(ii)\
If $\lambda= \rho-\mu$ then
 \begin{equation}\label{s6a}
 \int_{u\ge 0}
 E\left(e^{-\rho\Delta\Xbar_{\tau_u}-\mu\Xbar_{\tau_u-}
 -\nu \Gtau-\theta(\tau_u-\Gtau)}; \tau_u<\infty\right)\rmd u
=\frac{1}{\kappa(\nu,\mu)}
\frac{\partial_+ \kappa(\theta,\rho)}{\partial_+ \rho}
\end{equation}
provided the right derivative exists.
 \end{theorem}

The right derivative in \eqref{s6a} exists and equals the derivative if $\gk(\theta,\rho-\veps)<\infty$ for some $\veps>0$.
When $\lambda=0,\theta=\nu\ge 0, \rho>0, \mu>0$, \eqref{e3} reduces to
the {\it second factorization identity},
(3.2) of Percheskii and Rogozin (1969).
[Alili and Kyprianou (2005) give a short and elegant proof
of the second factorization identity using the strong Markov property.]
Theorem \ref{thmct3} can be used in the computation of certain exponential
Gerber-Shiu functionals from insurance risk, see Griffin and Maller (2010a).  Another application,  to stability of
the exit time, can be found in Griffin and Maller (2010b).

It is natural to ask if there is a
quintuple Laplace transform identity,  analogous to the quintuple law.
It is straightforward
to follow calculations similar to those in the proof of
Theorem \ref{thmct3}
and derive a corresponding expression to \eqref{e3},
but because the component
$X_{\tau_u-}$ cannot be expressed in terms of the ladder process,
the resulting expression cannot be expressed simply in terms of the
kappa functions.


\setcounter{equation}{0} \section{Bivariate Subordinators}\label{s4}

$(L^{-1},H)$ and $(\whL^{-1}, \whH)$ are, possibly killed,  bivariate subordinators,
and some of our results  require only this property.
In this  section we prove several theorems in this
generality.  These
will then be applied in Section \ref{s5}
 to the fluctuation variables.

In the fluctuation setting, let
\be\label{Tu}
T_u=T^{\cal H}_u=\inf \{s \geq0: {\cal H}_s>u\},\ u\ge 0.
\ee
Using $H_s=X_{L^{-1}_s}$ on $\{L^{-1}_s<\infty\}$, $s\ge 0$,
and recalling the exponential killing described in \eqref{H}, we have
\be\label{trans}
X_{\tau_u}= {\cal H}_{T_u},\
\Xbar_{\tau_u-}= {\cal H}_{T_u-},\
\tau_u= {\cal L}^{-1}_{T_u},\
 {\rm and}\
G_{\tau_u-}={\cal L}^{-1}_{T_u-}\
 {\rm on }\ \{T_u<e(q)\}.
\ee
Thus, via \eqref{ct2},
\be\label{pdef}
p(t,u)= P({\cal L}^{-1}_{T_u}\le t, {\cal H}_{T_u}=u, T_u<e(q)).
\ee
This suggests the following setup. Let $(Z,Y)$ be any two dimensional subordinator
obtained from a true subordinator
$({\cal Z},{\cal Y})$ by exponential killing at rate $q\ge 0$, say.
Corresponding  to \eqref{kapexp}, $(Z,Y)$ has Laplace exponent
$\tk(a,b)  =q-\log  E e^{-a{\cal Z}_1-b{\cal Y}_1}$ where
\be\label{tk}
\tk(a,b)=
q+\rmd_{Z}a+\rmd_Yb+\int_{t\ge0}\int_{x\ge0}
\left(1-e^{-at-bx}\right)
\Pi_{Z, Y}(\rmd t, \rmd x),
\ee
for values of $a,b\in \R$ for which the expression is finite.
Analogous to the  fluctuation variables, we define
\ben
T^{\cal Y}_u= \inf \{s \geq0: {\cal Y}_s>u\},\ u\ge 0,
\een
and
\be\label{defPZ}
p_{Z,Y}(t,u)= P({\cal Z}_{T^{\cal Y}_u}\le t, {\cal Y}_{T^{\cal Y}_u}=u, T^{\cal Y}_u<e(q)),
\ee
where $e(q)$ is an independent exponential random variable with parameter $q$.
Also set
\[
V_{Z,Y}(t,u)=\int_0^\infty e^{-qs}P({\cal Z}_s\le t,{\cal Y}_s\le u)\rmd s.
\]
So $V_{Z,Y}(\cdot,\cdot)$ has Laplace transform
\begin{equation}\label{LTVZY}
\int_{t\ge 0}\int_{x\ge 0}e^{-at-bx} V_{Z,Y}(\rmd t,\rmd x)
=\frac{1}{\kappa_{Z,Y}(a,b)}\ \text  { if } {\kappa_{Z,Y}(a,b)}>0.
 \end{equation}

 We begin by investigating  aspects of the regularity of  $p_{Z,Y}$ defined in \eqref{defPZ}.

 \begin{lemma}\label{sublemL_C2}\
The function $p_{Z,Y}(\cdot,\cdot)$ has  the following properties:

\begin{item} (a) $p_{Z,Y}(\cdot, u)$\text{ is  right continuous and
non-decreasing on  $[0,\infty)$ for every $u>0$};
\end{item}
\begin{item} (b) $p_{Z,Y}(t,\cdot)$\text{ is left continuous
 on $(0,\infty)$ for every } $t\ge 0$;
\end{item}
\begin{item} (c) \text{ If } $p_{Z,Y}(\cdot, u)$  \text{is
continuous on }  $(0,\infty)$ \text{for every } $u>0$,
\text{then }
$p_{Z,Y}(t,\cdot)$  \text{is  continuous on }  $(0,\infty)$
\text{ for every } $t> 0$.
\end{item}
 \end{lemma}

\bigskip \noindent {\bf Proof of Lemma \ref{sublemL_C2}:}\
First observe that the results trivially hold if
$\rmd_Y=0$, since then
${\cal Y}$ does not creep and so
$p_{Z,Y}(t,u)\le P({\cal Y}_{T^{\cal Y}_u}=u)=0$ for all $t\ge 0,u>0$.
Thus for the remainder of the proof we assume $\rmd_Y>0$.

Part (a) follows immediately from the definition of $p_{Z,Y}$.
To prove Parts (b) and (c) we will use the following two equations which
are simple consequences of the strong Markov property
(cf. Andrew (2006)): for any
$x>0,y>0,r\ge 0,s\ge 0$, we have
\begin{equation}\label{subap1}
p_{Z,Y}(r+s,x+y)\ge p_{Z,Y}(r,x)p_{Z,Y}(s,y)
\end{equation}
and
 \begin{equation}\label{subap2}
p_{Z,Y}(s,x+y)\le p_{Z,Y}(r,x)p_{Z,Y}(s,y)+1-p_{Z,Y}(r,x).
 \end{equation}

By Theorem III.5  of  Bertoin (1996), which applies to
nondefective subordinators with $\rmd_Y>0$,  we have
$\lim_{\varepsilon \downarrow 0}
P({\cal Y}_{T^{\cal Y}_\varepsilon}=\varepsilon)=1$.
Since $\rmd_Y>0$, $\cal {Y}$ is strictly increasing, and so  $T^{\cal Y}_\ve\downarrow 0$ a.s. as $\ve\dto 0$.
Thus for every $\delta>0$
\be\label{subpat0}
\lim_{\ve \downarrow 0}p_{Z,Y}(\delta,\ve)
=\lim_{\ve \downarrow 0}
P({\cal Z}_{T^{\cal Y}_\ve}\le \delta, {\cal Y}_{T^{\cal Y}_\ve}=\ve, T^{\cal Y}_\ve<e(q))=1.
\ee

Now fix $u>0$ and $t\ge 0$.  Then for any $0<\ve<u$ and $\delta>0$ we have by \eqref{subap1}  and \eqref{subap2}
\ben
{p_{Z,Y}(t,u)-1+p_{Z,Y}(\delta,\ve)}\le {p_{Z,Y}(\delta,\ve)}p_{Z,Y}(t,u-\ve)
\le  {p_{Z,Y}(t+\delta,u)}.
\een
Letting  $\ve\downarrow 0$ then $\delta\downarrow 0$, and using
Part (a) and \eqref{subpat0}, proves Part (b).
Similarly if in addition, $t>0$ and  $\delta<t$, then
\ben
p_{Z,Y}(\delta,\ve)p_{Z,Y}(t-\delta,u)\le p_{Z,Y}(t,u+\ve)\le
p_{Z,Y}(\delta,\ve)p_{Z,Y}(t,u)+ 1 -p_{Z,Y}(\delta,\ve).
\een
Letting $\ve\downarrow 0$ then $\delta\downarrow 0$, we conclude that
\ben
p_{Z,Y}(t-,u)\le \liminf_{\ve\downarrow 0} p_{Z,Y}(t,u+\ve)\le \limsup_{\ve\downarrow 0} p_{Z,Y}(t,u+\ve)\le p_{Z,Y}(t,u).
\een
Thus if
$p_{Z,Y}(\cdot,u)$ is continuous on $(0,\infty)$ for every $u>0$, then
$p_{Z,Y}(t,\cdot)$ is right continuous  on $(0,\infty)$ for every $t>0$.
Combining this with Part (b) proves $p_{Z,Y}(t,\cdot)$ is continuous
on $(0,\infty)$ for every $t>0$.  \hfill\halmos
\bigskip

\begin{lemma}\label{sublemL_C1}\
For any $u\ge0$ and $t\ge 0$,
\begin{equation}\label{subpint}
\int_0^up_{Z,Y}(t,v) \rmd v=\rmd_Y V_{Z,Y}(t,u).
 \end{equation}
 \end{lemma}

\bigskip \noindent {\bf Proof of Lemma \ref{sublemL_C1}:}\
If $d_Y=0$ then $\cal{Y}$ does not creep, so both sides of \eqref{subpint} are zero.  Thus we may assume $d_Y>0$.  First observe that for any $s$, $\{v:{\cal Y}_{T^{\cal Y}_v}=v, T^{\cal Y}_v=s\}$ is at most a singleton,
and so
\be\label{Yvv'}
\int_0^\infty{\bf 1}\{{\cal Y}_{T^{\cal Y}_v}=v, T^{\cal Y}_v=s\}\rmd v =0.
\ee
Next, if $T^{\cal Z}_t=\inf\{s:{\cal Z}_s>t\}$,
then
\be\label{ZS}
\{s<T^{\cal Z}_t\}\subset\{{\cal Z}_s\le t\}\subset\{s\le T^{\cal Z}_t\}.
\ee
Thus, using \eqref{Yvv'} and \eqref{ZS},
\ben\ba
\int_0^u p_{Z,Y}(t,v)\rmd v
&=\int_0^uP({\cal Z}_{T^{\cal Y}_v}\le t,
{\cal Y}_{T^{\cal Y}_v}=v, T^{\cal Y}_v<e(q))\rmd v\\
&= \int_0^uP({T^{\cal Y}_v}\le T^{\cal Z}_t, {\cal Y}_{T^{\cal Y}_v}=v,
T^{\cal Y}_v<e(q))\rmd v\\
&=\int_0^uP({\cal Y}_{T^{\cal Y}_v}=v,T^{\cal Y}_v\le T^{\cal Z}_t\wedge e(q))\rmd v.
\ea\een
Now since $d_Y>0$, $\cal{Y}$ is strictly increasing.  Thus if $\cal{Y}$ hits $v$ then it does so at time $T^{\cal Y}_v$.  Hence
\ben\ba
\int_0^u{\bf 1}\{{\cal Y}_{T^{\cal Y}_v}=v,T^{\cal Y}_v\le T^{\cal Z}_t\wedge e(q)\}\rmd v
={\cal Y}_{T^{\cal Y}_u\wedge T^{\cal Z}_t\wedge e(q)}-\sum_{s\le T^{\cal Y}_u\wedge
T^{\cal Z}_t\wedge e(q)}\Delta {\cal Y}_s
\ea\een
as each quantity represents the Lebesgue measure of the set of points in
$[0,u]$ hit by ${\cal Y}$ by time $T^{\cal Z}_t\wedge e(q)$.
Since
${\cal Y}_r=\rmd_Y r+ \sum_{s\le r}\Delta {\cal Y}_s$, this gives
 \be\label{pZY}
 \int_0^u p_{Z,Y}(t,v)\rmd v=\rmd_YE(T^{\cal Y}_u\wedge T^{\cal Z}_t\wedge e(q)).
 \ee
But by \eqref{ZS} (which also applies to $\cal {Y}$ and ${T^{\cal Y}_u}$)
\be\ba\label{ZY}
\int_0^\infty{\bf 1}\{{\cal Z}_{s}\le t,{\cal Y}_{s}\le u, s<e(q)\}\rmd s
&\le \int_0^\infty{\bf 1}\{s\le T^{\cal Z}_t,s\le T^{\cal Y}_u,   s<e(q)\}\rmd s\\
&=\int_0^\infty{\bf 1}\{s< T^{\cal Z}_t, s< T^{\cal Y}_u,  s<e(q)\}\rmd s\\
&\le \int_0^\infty{\bf 1}\{{\cal Z}_s\le t,{\cal Y}_s\le u, s<e(q)\}\rmd s.
 \ea\ee
Thus by \eqref{pZY} and \eqref{ZY},
\ben
\int_0^u p_{Z,Y}(t,v)\rmd v=
\rmd_Y\int_0^\infty P({\cal Z}_s\le t,{\cal Y}_s\le u, s<e(q))
\rmd s
=\rmd_YV_{Z,Y}(t,u).
\een
\hfill\halmos\bigskip

\begin{theorem}\label{subthmct1}\
Parts (i) and (ii) of Theorem \ref{thmct1}
hold precisely as stated with
$p_{Z,Y}$ in place of $p$,
$\rmd_Y$ in place of  $\rmd_H$,
and  $V_{Z,Y}$ in place of $V$.
\end{theorem}

\bigskip \noindent {\bf Proof of Theorem \ref{subthmct1}:}\
Since $V_{Z,Y}(t,u)>0$ for all $t>0$ and $u>0$ by right continuity of $(Z,Y)$, we have by
\eqref{subpint} and Lemma \ref{sublemL_C2} (b) that
$$
\rmd_Y>0
\text{ iff }\ \text{ $p_{Z,Y}(t,u)>0$ for some $t>0,u>0$}.
$$
On the other hand
\ben\ba
\rmd_Y>0
\text{ iff
$0<P({\cal Y}_{T^{\cal Y}_u}=u, T^{\cal Y}_u<e(q))=\lim_{t\to\infty}p_{Z,Y}(t,u)$ for every $u>0$}.
\ea\een
Combined with monotonicity of $p(\cdot,u)$ for $u>0$, these give
the equivalence of
the subordinator versions of
\eqref{C1}, \eqref{C2} and \eqref{C4}.
To complete the proof of Part (i) observe that
 the subordinator version of
\eqref{C4} implying  the subordinator version of
\eqref{C3} was proved in \eqref{subpat0},
while  the subordinator version of \eqref{C3} implying
 the subordinator version of \eqref{C1} is trivial.

If $d_Y>0$ then $\cal {Y}$ is not compound Poisson, so $V_{Z,Y}(t,0)=0$.  The remainder of part (ii) follows immediately from \eqref{subpint} and Lemma
\ref{sublemL_C2}(b).
\hfill\halmos
\bigskip

\begin{lemma}\label{subncreepPi}
For every $u>0$,
\be\label{subncreepPi1}
P(\Delta {\cal Z}_{T^{\cal Y}_u}>0, \Delta {\cal Y}_{T^{\cal Y}_u}=0)=0.
\ee
\end{lemma}

\bigskip \noindent {\bf Proof of Lemma \ref{subncreepPi}:}\
If ${\cal Y}$ is compound Poisson,
then $P(\Delta {\cal Y}_{T^{\cal Y}_u}=0)=0$, so the  result is trivial.
Thus we may assume that ${\cal Y}$ is not compound Poisson, in which case ${\cal Y}$ is
strictly increasing.
Thus by the compensation formula, p7 of Bertoin (1996),
for every $\ve>0$
\ben\ba
P(\Delta {\cal Y}_{T^{\cal Y}_u}=0, \Delta {\cal Z}_{T^{\cal Y}_u}>\ve)
&=E\sum_{t>0}{\bf 1}\{{\cal Y}_{t-}=u, \Delta {\cal Y}_{t}=0, \Delta {\cal Z}_{t}>\ve\}
\\
&=\Pi_{Z,Y}((\ve,\infty)\times\{0\})\int_0^\infty P({\cal Y}_t=u)\rmd t.
\ea\een
The last expression is 0
by
Proposition I.15 of Bertoin (1996).\hfill\halmos\bigskip

The next result is a ``quadruple law" for $(Z,Y)$,
in a similar spirit to the quintuple law.

\begin{theorem}[Quadruple Law] \label{quadlaw}\
For $u>0$,  $ x\ge 0$, $0\le y\leq u$, $s\ge 0$, $t\ge 0$, we have
\begin{eqnarray}\label{subDKc}
&&P\left({\cal Y}_{T^{\cal Y}_u}-u \in \rmd x,
u-{\cal Y}_{T^{\cal Y}_u-} \in \rmd y,
\Delta {\cal Z}_{T^{\cal Y}_u} \in \rmd s, {\cal Z}_{T^{\cal Y}_u-} \in \rmd t; T^{\cal Y}_u<e(q)\right)
\nonumber\\
&&
={\bf 1}_{\{x>0\}}
|V_{Z,Y}(\rmd t,u-\rmd y)|\Pi_{Z,Y}(\rmd s, \rmd x+y)
+\rmd_Y \frac{\partial_-}{\partial_- u}V_{Z,Y}(\rmd t,u)
\delta_0(\rmd s,\rmd  x,\rmd y),
\nonumber\\
\end{eqnarray}
with the convention that the term containing the differential
${\partial_-}V_{Z,Y}(\rmd t,u)/{\partial_- u}$ is absent when
$ \rmd_Y=0$
(in which case  ${\partial_-}V_{Z,Y}(t,u)/{\partial_- u}$
need not be defined).
\end{theorem}

\bigskip \noindent {\bf Proof of Theorem \ref{quadlaw}:}\
Fix $u>0$. By the compensation formula,
we get for $ x>0$, $0\le y\leq u$, $s\ge 0$, $t\ge 0$,
\be\ba\label{subDK}
P({\cal Y}_{T^{\cal Y}_u}-u &\in  \rmd x,u-{\cal Y}_{T^{\cal Y}_u-} \in \rmd y,
\Delta {\cal Z}_{T^{\cal Y}_u} \in \rmd s, {\cal Z}_{T^{\cal Y}_u-} \in \rmd t; T^{\cal Y}_u<e(q))\\
&=P(\Delta {\cal Y}_{T^{\cal Y}_u} \in \rmd x +y,u-{\cal Y}_{T^{\cal Y}_u-} \in \rmd y,
\Delta {\cal Z}_{T^{\cal Y}_u} \in \rmd s, {\cal Z}_{T^{\cal Y}_u-} \in \rmd t; T^{\cal Y}_u<e(q))\\
&=E\sum_{r>0} {\bf 1}\{\Delta {\cal Y}_{r} \in \rmd x +y,u-{\cal Y}_{r-} \in \rmd y,
\Delta {\cal Z}_{r} \in \rmd s, {\cal Z}_{r-} \in \rmd t, r<e(q)\}\\
&=E\int_0^\infty {\bf 1}\{u-{\cal Y}_{r-} \in \rmd y, {\cal Z}_{r-} \in \rmd t,r<e(q)\}
\rmd r\, \Pi_{Z,H}(\rmd s, \rmd x +y)\\
&=|V_{Z,H}(\rmd t,u-\rmd y)|\Pi_{Z,H}(\rmd s, \rmd x+y).
\ea\ee
Thus we  have left to consider the case $x=0$.

First observe that from  Proposition III.2 of Bertoin (1996), it follows that
\begin{equation}\label{suballu-}
P({\cal Y}_{T^{\cal Y}_u-}<u={\cal Y}_{T^{\cal Y}_u}) =0, \ {\rm for\ all}\ u>0.
 \end{equation}
Now suppose  $\rmd_Y>0$.
By part (ii) of Theorem \ref{subthmct1}
\begin{equation}\label{subct1}
p_{Z,Y}(t,u)
=\rmd_Y \frac{\partial_-}{\partial_- u}V_{Z,Y}(t,u).
\end{equation}
This, together with \eqref{subncreepPi1} and \eqref{suballu-}, shows
\be\ba\label{subct1a}
P({\cal Y}_{T^{\cal Y}_u}=u, &u-{\cal Y}_{T^{\cal Y}_u-} \in \rmd y,
\Delta {\cal Z}_{T^{\cal Y}_u}\in \rmd s,  {\cal Z}_{T^{\cal Y}_u-} \in \rmd t; T^{\cal Y}_u<e(q) )\\
&=
P({\cal Y}_{T^{\cal Y}_u}=u,   {\cal Z}_{T^{\cal Y}_u} \in \rmd t; T^{\cal Y}_u<e(q)  )\delta_0(\rmd s,\rmd y)\\
&=
\rmd_Y \frac{\partial_-}{\partial_- u}V_{Z,Y}(\rmd t,u)
\delta_0(\rmd s,\rmd y)
\ea\ee
for  $0\le y\leq u$, $s \geq0$ and $t\ge 0$.
When $\rmd_Y=0$, \eqref{subct1a} continues to hold since the  lefthand side
of
\eqref{subct1a} is 0 because ${\cal Y}$ does not creep.
Adding \eqref{subct1a} to \eqref{subDK} then  gives \eqref{subDKc}.
\newline{}
\hfill\halmos\bigskip

The next result is a generalisation of
Theorem \ref{thmct3} to the subordinator setup.  The conditions on $\tk$ are analogous to those on $\gk$ in Theorem \ref{thmct3}.

\begin{theorem}[A Laplace Transform Identity for Subordinators]
\label{subthmct3}\
\newline
Fix $\mu, \rho, \lambda, \nu, \theta$ so that $\tk(\theta,\mu+\lambda), \tk(\theta,\rho)$ are finite and $\tk(\nu,\mu)> 0$.

\noindent
(i)\
If $\lambda\ne \rho-\mu$ then
\begin{eqnarray}\label{sube3}
&&\int_{u\ge 0}e^{-\mu u}
E\left(e^{-\rho({\cal Y}_{T^{\cal Y}_u}-u)-\lambda(u-{\cal Y}_{T^{\cal Y}_u-})
-\nu {\cal Z}_{T^{\cal Y}_u-}-\theta\Delta {\cal Z}_{T^{\cal Y}_u}}; T^{\cal Y}_u<e(q) \right)\rmd u
\nonumber\\
&&\qquad \qquad \qquad \qquad \qquad \qquad
=
\frac {\tk(\theta,\mu+\lambda)-\tk(\theta,\rho)}
{(\mu+\lambda-\rho)\tk(\nu,\mu)}.
\end{eqnarray}
(ii)\
If $\lambda=\rho-\mu$ then
 \begin{equation}\label{subs6a}
 \int_{u\ge 0}
 E\left(e^{-\rho\Delta {\cal Y}_{T^{\cal Y}_u}-\mu {\cal Y}_{T^{\cal Y}_u-}
 -\nu {\cal Z}_{T^{\cal Y}_u-}-\theta \Delta {\cal Z}_{T^{\cal Y}_u}}; T^{\cal Y}_u<e(q) \right)\rmd u
=\frac{1}{\kappa_{Z,Y}(\nu,\mu)}
\frac{\partial_+\kappa_{Z,Y}(\theta,\rho)}{\partial_+ \rho}
\end{equation}
provided the right derivative exists.
 \end{theorem}

\bigskip \noindent {\bf Proof of Theorem \ref{subthmct3}:}\
{}Taking the expectation over the set $\{{\cal Y}_{T^{\cal Y}_u}>u\}$, from \eqref{subDKc} we find

\be\ba\label{sube1}
&E\left(e^{-\rho({\cal Y}_{T^{\cal Y}_u}-u)-\lambda(u-{\cal Y}_{T^{\cal Y}_u-})
-\nu {\cal Z}_{T^{\cal Y}_u-}-\theta(\Delta {\cal Z}_{T^{\cal Y}_u})}; T^{\cal Y}_u<e(q),
{\cal Y}_{T^{\cal Y}_u}>u\right) \\
&\qquad =
\int_{0\le y\le u}\int_{ x>0}\int_{s\ge 0}\int_{t\ge 0}
e^{-\rho x-\lambda y-\nu t-\theta s}|V_{Z,Y}(\rmd t,u-\rmd y)|
\Pi_{Z,Y}(\rmd s, \rmd x+y)\\
&\qquad =
\int_{0\le w\le u}\int_{x>u-w}\int_{s\ge 0}\int_{t\ge 0}
e^{-\rho(x-u+w)-\lambda (u-w)-\nu t-\theta s}V_{Z,Y}(\rmd t,\rmd w)
\Pi_{Z,Y}(\rmd s, \rmd x).
\ea\ee
Now take the Laplace transform of both sides of
\eqref{sube1}. For $\lambda\ne \rho-\mu$, we obtain
\be\ba\label{sube2a}
\int_{u\ge 0}e^{-\mu u}\rmd u\int_{0\le w\le u}
&\int_{x>u-w}e^{-\rho(x-u+w)-\lambda (u-w)}\\
&\qquad=
\int_{w\ge 0}\int_{x>0}\int_{w\le u< w+x}
e^{-(\mu+\lambda-\rho)u}e^{-(\rho-\lambda)w}e^{-\rho x}\rmd u\\
&\qquad=
\frac 1{\rho-\mu-\lambda}\int_{w\ge 0}
e^{-\mu w}\int_{x> 0}(e^{-(\mu+\lambda)x}-e^{-\rho x}).
\ea\ee
Since we may clearly also include $x=0$ in the last integral,  we then have
\be\ba\label{sube3a}
&\int_{u\ge 0}e^{-\mu u}
E\left(e^{-\rho({\cal Y}_{T^{\cal Y}_u}-u)-\lambda(u-{\cal Y}_{T^{\cal Y}_u-})
-\nu {\cal Z}_{T^{\cal Y}_u-}-\theta(\Delta {\cal Z}_{T^{\cal Y}_u})}; T^{\cal Y}_u<e(q),
{\cal Y}_{T^{\cal Y}_u}>u\right)\rmd u
\\
&\qquad=
\frac 1{\rho-\mu-\lambda}\int_{w\ge 0}\int_{t\ge 0}
e^{-\nu t-\mu w}V_{Z,Y}(\rmd t, \rmd w)
\\
&\qquad\qquad\qquad \qquad\qquad\qquad \qquad\qquad
\times
\int_{s\ge 0}\int_{x\ge 0}
e^{-\theta s}(e^{-(\mu+\lambda)x}-e^{-\rho x})
\Pi_{Z,Y}(\rmd s, \rmd x)\\
&\qquad=\frac {\tk(\theta,\mu+\lambda)-\tk(\theta,\rho)-
(\mu+\lambda-\rho)\rmd_Y}{(\mu+\lambda-\rho)\tk(\nu,\mu)}\\
&\qquad=
\frac {\tk(\theta,\mu+\lambda)-\tk(\theta,\rho)}
{(\mu+\lambda-\rho)
\tk(\nu,\mu)}-\frac {\rmd_Y}{\tk(\nu,\mu)}
\ea\ee
by \eqref{tk} and \eqref{LTVZY}.

If $\rmd_Y>0$ then we need to add in the second term in \eqref{subDKc} due to creeping.  From \eqref{subDKc} and part (ii) of Theorem \ref{subthmct1}
we have
\ben\ba
\int_{u\ge 0}e^{-\mu u}&
E\left(e^{-\rho({\cal Y}_{T^{\cal Y}_u}-u)-\lambda(u-{\cal Y}_{T^{\cal Y}_u-})
-\nu {\cal Z}_{T^{\cal Y}_u-}-\theta(\Delta {\cal Z}_{T^{\cal Y}_u})}; T^{\cal Y}_u<e(q), {\cal Y}_{T^{\cal Y}_u}=u\right)\rmd u
\\
&\qquad
=\rmd_Y \int_{u\ge 0}e^{-\mu u}\int_{t\ge 0}
e^{-\nu t}\frac{\partial_-}{\partial_- u}V_{Z,Y}(\rmd t,u)\rmd u
\\
&\qquad
=\rmd_Y\int_{u\ge 0}e^{-\mu u}\int_{t\ge 0}
e^{-\nu t}V_{Z,Y}(\rmd t,\rmd u)
\\
&\qquad=
\frac {\rmd_Y}{\tk(\nu,\mu)}.
\ea\een
Added to \eqref{sube3a}, this gives  \eqref{sube3}.

Now consider the case where $\lambda=\rho-\mu$.
Let $\ve>0$ and set  $\lambda'=\rho-\mu+\ve$.  Then $\tk(\theta,\mu+\lambda')$ is finite since $\tk(\theta,\rho )$
is finite. Thus by \eqref{sube3}
 \begin{equation*}\ba
\int_{u\ge 0}
&E\left(e^{-\rho\Delta {\cal Y}_{T^{\cal Y}_u}-\mu {\cal Y}_{T^{\cal Y}_u-}
-\ve(u-{\cal Y}_{T^{\cal Y}_u-})
 -\nu {\cal Z}_{T^{\cal Y}_u}-\theta\Delta  {\cal Z}_{T^{\cal Y}_u}}; T^{\cal Y}_u<e(q)\right)\rmd u\\
&\qquad\qquad\qquad\qquad\qquad\qquad\qquad=\frac {\kappa_{Z,Y}(\theta,\rho+\ve)-\kappa_{Z,Y}(\theta,\rho)}
{\ve\kappa_{Z,Y}(\nu,\mu)}.
\ea\end{equation*}
Letting $\ve\downarrow0$ and using monotone convergence completes the
proof of \eqref{subs6a}.
\hfill\halmos\bigskip

Results similar to \eqref{sube3} can be found in Winkel (2005).  
Winkel's interest in bivariate subordinators is in modelling electronic foreign exchange markets and he does not make the connection with the ladder height process.
As mentioned in the introduction, it appears to have gone previously  unnoticed, that the second factorization identity is a special case of a general transform result for bivariate subordinators.


\setcounter{equation}{0}
\section{Proofs for Section \ref{s3}}
\label{s5}
We now turn to the proofs of the fluctuation results from
Section \ref{s3}.
Recall  the definitions of $p(t,u)$ and $T_u$ in \eqref{ct2}
and \eqref{Tu} respectively, and
from \eqref{pdef} that
\ben
p(t,u)= P({\cal L}^{-1}_{T_u}\le t, {\cal H}_{T_u}=u, T_u<e(q)).
\een
In view of the correspondences
$ (L^{-1},H)\leftrightarrow (Z,Y)$,
$ p(t,u)\leftrightarrow p_{Z,Y}(t,u)$,
and
$ V(t,u)\leftrightarrow V_{Z,Y}(t,u)$, we can carry a number of results directly
across from Section \ref{s4}.
In particular, the results of
Lemma \ref{sublemL_C2} and Lemma \ref{sublemL_C1}
hold with  $p$ in place of $p_{Z,Y}$
and  $V$ in place of $V_{Z,Y}$.

\bigskip \noindent {\bf Proof of Theorem \ref{thmct1}:}\
Parts (i) and (ii) follow
immediately from Theorem  \ref{subthmct1}.
For Part (iii), supposing
$X$ is not compound Poisson with positive drift,
by Theorem 27.4 of Sato (1999),
\ben
P(\tau_u= t, X_{\tau_u}=u)\le P(X_{t}=u)=0
\een
for all $u>0$.  Consequently $p(\cdot, u)$ is  continuous on
$[0,\infty)$ for every $u>0$, and hence by Lemma \ref{sublemL_C2},
$p(t, \cdot)$ is  continuous on $(0,\infty)$ for every $t> 0$.  Since
$p(0, \cdot)\equiv 0$ on $(0,\infty)$ this continues to hold for $t=0$.
By \eqref{subpint},
this then implies differentiability of $V(t,\cdot)$.
\hfill\halmos
\bigskip

\begin{example}\label{E2}
{\rm Let $X_t=t+Y_t$ where $Y$ is compound Poisson with L\'evy
measure $\Pi_Y(\rmd x)
= \delta_{\{1\}}(\rmd x) + \delta_{\{-1\}}(\rmd x)$.  Since
$X_t=t+\sum_{s\le t}\Delta Y_s$ and $\sum_{s\le t}\Delta Y_s$
is an integer, we have that for any $t\in (0,1)$, $u>0$
\ben
p(t,u)>0  \text{ iff } u\in \cup_{n=0}^\infty (n, n+t].
\een
Furthermore $p(t,u)\ge e^{-2t u}$ for $u\in (0,t]$.
Thus neither $p(t,\cdot)$ nor $p(\cdot,u)$ is  continuous,
and, using \eqref{subpint},
$V(t,\cdot)$ is not differentiable.  Finally,
in contrast to  monotonicity of $p(\cdot,u)$,
  $p(t,\cdot)$ is not monotone.
}\end{example}

\begin{example}\label{E3} {\rm
Set $X_t=t-Y_t$,
 where $Y$ is a subordinator and $\Pi_Y(\R)=\infty$.
Then clearly $p(t,u)=0$ for $u>t$.
Thus there is no hope of proving $p(t,u)>0$ for all $t,u>0$
even in the situation of (iii) of Theorem \ref{thmct1}.
}\end{example}
\bigskip

Before turning to the proof of Theorem \ref{thmct2},
we need a preliminary result, generalising \eqref{suballu-}, which is surely well known, but for which we can not find an exact reference.

\begin{lemma}\label{nj}  Fix $u>0$.

\noindent
(i)\ If  $X$ is not compound Poisson, then
\begin{equation}\label{not}
P(X_{t-}\neq u, X_{t}=u\ {\rm for\ some }\ t>0)
=P(X_{t-}=u, X_{t}\neq u\ {\rm for\ some }\ t>0)=0;
\end{equation}
(ii)\ For any L\'evy process $X$ and $u>0$
\begin{equation}\label{allu-}
P(X_{\tau_u-} <u,X_{\tau_u}=u, \tau_u<\infty)=0.
\end{equation}
\end{lemma}


\bigskip \noindent {\bf Proof of Lemma \ref{nj}:}\
(i)\
Use the compensation formula  to write
for $u>0$
\ben\ba
P(X_{t-}\neq u, X_{t}=u\ {\rm for\ some }\ t>0)
&\le
E\sum_{t>0}{\bf 1}\{X_{t-}\neq u, X_{t-}+\Delta X_{t}=u\}\\
&=
E\int_0^\infty \rmd t
\int_{\xi\ne 0} {\bf 1}\{X_{t-}\neq u, X_{t-}+\xi=u\}\Pi_X(d\xi)\\
&=
\int_{\xi\ne 0}\int_0^\infty  P(X_{t-}=u-\xi)\rmd t\,
\Pi_X(d\xi).
\ea\een
The last expression is 0 when $X$ is not compound Poisson,
since the potential measure of $X$ is diffuse by Proposition I.15 of Bertoin (1996).  Similarly for any
$\ve>0$
\ben\ba
P(X_{t-}=u, |X_{t}-u|>\ve\ {\rm for\ some }\ t>0)
&\le
E\sum_{t>0}{\bf 1}\{X_{t-}= u, |\Delta X_{t}|>\ve\}\\
&=
\Pi_X([-\ve,\ve]^c)\int_0^\infty
P(X_{t-}= u)\ \rmd t=0.
\ea\een
Letting $\ve\to 0$ completes the proof.

(ii)\
Clearly we may assume $P(X_{\tau_u}= u)>0$ else there is nothing to
prove. In that case $X$ is not compound Poisson. Since
\ben
\{X_{\tau_u-}<u, X_{\tau_u}=u, \tau_u<\infty\}\subseteq
\{X_{t-}\neq u, X_{t}=u\ {\rm for\ some }\ t>0\},
\een
\eqref{allu-} follows from \eqref{not}.
\hfill\halmos
\bigskip

\noindent {\bf Proof of Theorem \ref{thmct2}:}\
Fix  $v\ge 0$, $0\le y\leq u \wedge v$, $s \geq0$ and $t\ge 0$.  For $x>0$ this is Doney and Kyprianou's quintuple law.
If $x=0$,  then from \eqref{trans}, \eqref{subDKc} and \eqref{allu-}
\ben\ba
&P\left(X_{\tau_u}-u\in \rmd x, u-X_{\tau_u-} \in \rmd v,
u-\overline{X}_{\tau_u-}  \in \rmd y,
\tau_u-\Gtau  \in \rmd s, \Gtau  \in \rmd t \right)\\
&\qquad  =
P\left({\cal H}_{T_u}-u\in  \rmd x, u-{\cal H}_{T_u-}  \in \rmd y, \Delta {\cal L}^{-1}_{T_u}\in \rmd s, {\cal L}^{-1}_{{T_u}-}\in  \rmd t; T_u<e(q) \right)
\delta_0(\rmd v)\\
&\qquad  =
\rmd_H \frac{\partial_-}{\partial_- u}V(\rmd t,u)
\delta_0(\rmd s,\rmd x,\rmd v,\rmd y).
\ea\een
\hfill\halmos

\bigskip \noindent {\bf Proof of Corollary \ref{jtop}:}\
This follows easily from Theorem \ref{thmct2} and  \eqref{DKcor6}.  Alternatively use Theorem \ref{quadlaw}.
\hfill\halmos

\bigskip \noindent {\bf Proof of Theorem \ref{creepPi}:}\
Assume $X$ is not compound Poisson.
We can first easily dispense with the case that $0$ is irregular for
$(0,\infty)$; because then, $0$ is irregular for
$[0,\infty)$, so  by construction,
$L^{-1}$ and $H$  jump at the same times  (see p.24 of Doney (2005)).
Hence
$\Pi_{L^{-1},H}(\rmd s, \{0\})$ is the zero measure.
On the other hand,
by the strong Markov property at time
$\tau_u$,  $X$ does not creep over any $u\ge 0$.
Thus the result holds in this case.

We now assume that $0$ is regular for $(0,\infty)$, in which case the closure of the zero set of $\Xbar-X$ is a
perfect nowhere dense set with probability one; see for example the discussion on p.104 of Bertoin (1996).
Fix $r\in \Q, r>0$ and set $Y^r_t=X_{r+t}-X_r$. Let
\ben
\tau^r=\inf\{t>0:Y^r_t>\Xbar_{r}-X_r\}
\een
and
\ben A_r=\{Y^r_{\tau^r}=\Xbar_{r}-X_r,\ \Xbar_r-X_r>0\}.\een
By independence
\ben
P(A_r)=\int_{u> 0}P(X_{\tau_u} = u)P(\Xbar_r-X_r\in \rm du).
\een
Thus $P(A_r)>0$ iff $X$ creeps.
On the other hand, $\Pi_{L^{-1},H}(\rmd s, \{0\})$ is the zero measure
iff
$P( \Delta L^{-1}_s>0, \Delta H_s=0\text { for some } s>0)=0$.
The result then follows from the key observation that
\ben
\{\Delta L^{-1}_s>0, \Delta H_s=0
\text { for some } s>0\} \overset{a.s.}= \bigcup_{r\in\Q}A_r.
\een
To see this, assume $s$ is such that
$\Delta H_s=0$ and  $\Delta L^{-1}_s>0$.
Let $r\in (L^{-1}_{s-},L^{-1}_s)\cap \Q$.
Then $ \Xbar_r-X_r>0$, $X_{L^{-1}_s}= \Xbar_r$ and $X_t< \Xbar_r$ for
$L^{-1}_{s-}<t<L^{-1}_s$.
Since  a.s. the zero set of $\Xbar-X$ contains no
isolated points, it follows that off this $P$-null set
$Y^r_{\tau^r}=\Xbar_r-X_r$.
Thus $A_r$ occurs.
Conversely fix $r\in\Q$ and assume that
$Y^r_{\tau^r}=\Xbar_r-X_r>0$.  Choose $s>0$ so that
$r\in (L^{-1}_{s-},L^{-1}_s)$.
Then on $A_r$, $H_s=X_{L^{-1}_s}= \Xbar_r=H_{s-}$.
Hence $\Delta H_s=0$ and $\Delta L^{-1}_s>0$.
\hfill\halmos
\bigskip

As remarked earlier, the jumps of $(L^{-1},H)$ for which
$\Delta H=0$, do not occur when $X$ creeps over a fixed level.
This follows from an application of Lemma \ref{subncreepPi},
which gives
\be\label{ncreepPi1}
P(\Delta {\cal L}^{-1}_{T_u}>0, \Delta {\cal H}_{T_u}=0)=0, \ u>0.
\ee

Let
\ben
V_{H}(\rmd v)= \int_{t\ge 0}V(\rmd t, \rmd v)=\int_0^\infty P(H_s\in \rmd v)\rmd s
\een
be the potential measure of $H$, and similarly for $\whV_{\whH}$.

\bigskip \noindent {\bf Proof of Theorem \ref{CDK}}\
Assume $X$ is not compound Poisson.  Then $\whV_{\whH}$ is diffuse on $(0,\infty)$ by
Lemma 1 of Chaumont and Doney (2010). Then since $\Pi_X(\{0\})=0$ we have
\be\label{qq}
\int_{s\ge 0}\int_{v\ge 0}\whV(\rmd s, \rmd v)\Pi_X(\{v\})=\int_{v\ge 0}\whV_{\whH}(\rmd v)\Pi_X(\{v\})=0.
\ee
Now by \eqref{DKcor6},
\be
\Pi_{L^{-1},H}(\rmd s,\rmd x)
=\int_{v\ge 0} \whV(\rmd s, \rmd v)\Pi_X(\rmd x+v) \ {\rm{for\ all}}\ s\ge 0, x\ge 0
\ee
is equivalent to
\be
\Pi_{L^{-1},H}(\rmd s,\{0\})
=\int_{v\ge 0} \whV(\rmd s, \rmd v)\Pi_X(\{v\}) \ {\rm{for\ all}}\ s\ge 0
\ee
which in turn, by \eqref{qq}, is equivalent to $\Pi_{L^{-1},H}(\rmd s,\{0\})$ being the zero measure.
This is equivalent to $X$  not creeping by Theorem \ref{creepPi}.

If $X$ is compound Poisson, then recalling the definition of $(\whL^{-1}, \whH)$ prior to \eqref{WH} in this case, a natural approach is to apply the result
in the non compound Poisson case to the approximating process $X^\ve_t=X_t+\varepsilon t$ and take the limit as $\varepsilon\downarrow 0$.  Unfortunately this
cannot work since $X^\ve$ creeps and so the non compound Poisson result does not apply.
Consequently we are forced to appeal to a direct construction of $(\whL^{-1}, \whH)$ in this case.  We defer the details of this to the appendix.
\hfill\halmos

\begin{example} {\rm Assume that $X$ is a spectrally
negative compound Poisson process with positive drift $d_X$. For simplicity,
to avoid killing, also assume $EX_1\ge 0$. In this example $X$ creeps, and it is a simple matter to find $\Pi_{L^{-1},H}(\rmd s, \{0\})$.
Choose the normalisation of local time so that
\ben
L_t= \int_0^t {\bf 1}(X_r=\Xbar_r)\rmd r.
\een
Let $\sigma_1$ be the time of the first jump of $X$ and
\be\label{rho}
\alpha=\inf\{t>\sigma_1:
X_t\ge \Xbar_{\sigma_1}\}.
\ee
Then it is easy to see that (draw a picture)
\ben
L^{-1}_s=s+\sum_{i=1}^{N_s}R_i, \quad H_s=d_X s,
\een
where  $N$ is a Poisson process of rate $\lambda=\Pi_X(\R)$
independent of the i.i.d.
sequence $R_i, i\ge 1$, where
$R_1\overset{d}= \alpha-\sigma_1$.  From this we conclude
\ben
\Pi_{L^{-1},H}(\rmd s, \{0\})=
\lambda P(R_1\in \rmd s)=\Pi_{L^{-1}}(\rmd s), s>0.
\een
On the other hand, since $X$ is not compound Poisson, \eqref{qq} shows
that the lefthand side of \eqref{DKcor6} is the zero measure when $x=0$.

}\end{example}
\bigskip

\bigskip \noindent {\bf Proof of Theorem \ref{thmct3}:}\
Upon using \eqref{trans}, this is an immediate application of Theorem \ref{subthmct3}.  \hfill\halmos\bigskip


\setcounter{equation}{0} \section{Appendix; Proof of Theorem \ref{CDK} in the
Compound Poisson Case}\label{s6}

We want to prove Theorem \ref{CDK} in the compound Poisson case.  As mentioned earlier the natural approach using the approximating process
$X^\ve_t=X_t+\varepsilon t$ does not work.  Thus we proceed by using the random walk embedded in $X$ to give a direct construction of the
bivariate ladder processes.  

Throughout this section we assume that  $X$ is  compound Poisson. We write  it in the form
$X_t=\sum_{k=1}^{N_t}Y_k$,
where $N$ is a Poisson process of rate $\lambda>0$, $Y_k$ are
i.i.d. rvs  independent of $N$ with distribution function
$F$, and $P(Y_k=0)=0$.
Thus the L\'evy measure  of $X$ is $\Pi_X(\rmd x)=\lambda F(\rmd x)$.

Let $S_0=0$, and
$S_n=\sum_{k=1}^nY_k$,
$n\ge 1$.
Let $\sigma_0=0$ and $\sigma_k$, $k\ge 1$, denote the successive
jump times  of $X_t$.  Then
$\{S_n, n\ge 0\}$ is independent of $\{\sigma_k, k\ge 0\}$ and
$X_t=S_n$ for $\sigma_n\le t< \sigma_{n+1}$, $n\ge 0$.
Let $t_0=0$,
\be\label{walt}
t_{n+1}=\min\{m>t_n: S_m\ge S_{t_n}\},\quad n=0,1,2, \ldots,
\ee
be the {\em weak ascending} ladder times of $S_n$,
$h_n=S_{t_n}$ the corresponding ladder height sequence, and
\ben
U(k, x)=\sum_{n\ge 0}P(t_n=k, h_n\le x), \ k\ge 0, x\ge 0,
\een
the  corresponding bivariate renewal measure. The {\em strict ascending} ladder process is obtained by replacing the
inequality in \eqref{walt} with a strict inequality.
The analogous quantities for the {\em strict  descending}
ladder process will be denoted with a hat;  thus
\ben
\whU(k, x)=\sum_{n\ge 0}P(\wht_n=k, \whh_n\le x), \
k\ge 0, x\ge 0.
\een

We choose  the normalisation of the local time $L$ so that
\be\ba\label{LT}
L_t&= \int_0^t {\bf 1}(X_r=\Xbar_r)\rmd r=  \int_0^t \sum_{n=0}^\infty
{\bf 1}_{[\sigma_{t_n}, \sigma_{t_n+1})}(r)\rmd r, \ t\ge 0.
 \ea\ee
As remarked earlier this gives rise to the
{\it weak bivariate ladder process} $(L^{-1}, H)$ of $X$.
For $\whX$, we require $(\whL^{-1}, \whH)$ to be the
{\it strict} bivariate ladder process,
which  may be viewed as the limit of the  ascending bivariate ladder
process of $\whX_t-\ve t$ as $\ve\downarrow 0$.
For a  direct construction of $(\whL^{-1}, \whH)$, let $M_s$ be an auxiliary
Poisson process of rate 1, independent of $X$, and set
\be\label{wlp}
(\whL_s^{-1}, \whH_s)=\begin{cases} (\sigma_{\wht_n}, \whh_n) & \text{ if } M_s=n, \wht_n<\infty, n\ge 0,\\
(\infty,\infty) & \text{ if } M_s=n, \wht_n=\infty, n\ge 0.
\end{cases}
\ee
This is analogous to the definition of the ascending ladder processes
of $\whX$ in the non-compound Poisson case
when 0 is irregular for $[0,\infty)$  for $\whX$; see for example p.24 of Doney (2005).
To emphasize, $\whV(\cdot,\cdot)$ and $\widehat\kappa(\cdot,\cdot)$
are defined in terms of $(\whL^{-1}, \whH)$ given by \eqref{wlp}.
We should also remark that with these normalisations of the local times,
one can check that \eqref{WH} holds.

The connection between the bivariate ladder processes of $X$ and $S$
is given by

\begin{lemma}\label{lemVU}  For any $t\ge 0, x\ge 0$,
\be\ba\label{VU}
V(t,x)&=\lambda^{-1}\sum_{k\ge 0}
P(\sigma_{k+1}\le t)U(k, x),
\ea\ee
and
\be\label{tV0}
\whV(t, x)
=\sum_{k\ge 0}P( \sigma_{k}\le t)\widehat U(k,x).
\ee
\end{lemma}

\bigskip \noindent {\bf Proof of Lemma \ref{lemVU}:}\
Since \eqref{VU} will not be used in the sequel, we only prove  \eqref{tV0}.
We have
\be\ba\label{whVU1}
 P(\whL^{-1}_s\le t, \whH_s\le x)
&=\sum_{n=0}^\infty P(\sigma_{\wht_n}\le t, \whh_n\le x, M_s=n)\\
&=\sum_{n=0}^\infty \sum_{k=0}^\infty P(\wht_n=k, \whh_n\le x,
\sigma_{k}\le t, M_s=n)\\
&=\sum_{k=0}^\infty \sum_{n=0}^\infty
P(\wht_n=k, \whh_n\le x)P(\sigma_{k}\le t)P( M_s=n).
\ea\ee
One easily checks that $\int_0^\infty P( M_s=n)\ \rmd s =1$
for every $n\ge 0$, hence integrating \eqref{whVU1} we obtain
\ben\ba
\whV( t,  x)
&=\sum_{k=0}^\infty \sum_{n=0}^\infty
P(\wht_n=k, \whh_n\le x)P(\sigma_{k}\le t)\\
&=\sum_{k\ge 0}P( \sigma_{k}\le t)\widehat U(k,x),
\ea\een
as required.
\hfill\halmos
\bigskip

Recall the definition of $\alpha$ in \eqref{rho},
which in the present case reduces to
\ben
\alpha=\inf\{t>\sigma_1: X_t\ge 0\}.
\een
Also introduce
\ben
\beta=\inf\{n>0: S_n\ge 0\}.
\een

\begin{lemma}\label{tVPX}
For $s> 0, x\ge 0$ and $v\ge 0$
\be\label{XtV}
\whV(\rmd s, \rmd v)F(\rmd x+v)=P(X_{\ga-}\in -\rmd v, X_{\ga}\in \rmd x, \ga-\sigma_1\in \rmd s).
\ee
\end{lemma}

\bigskip \noindent {\bf Proof of Lemma \ref{tVPX}:}\
Using $P(\sigma_0\in \rmd s)=0$
when $s>0$ in the fourth equality below, duality in the fifth, and \eqref{tV0} in the sixth, we have
\ben\ba
P(X_{\ga-}\in -\rmd v, &X_{\ga}\in \rmd x, \ga-\sigma_1\in \rmd s)\\
&={\bf 1}_{\{v>0\}}P(X_{\ga-}\in -\rmd v, X_{\ga}\in \rmd x, \ga-\sigma_1\in \rmd s)\\
&={\bf 1}_{\{v>0\}}\sum_{i\ge 2}P(\gb=i, S_{i-1}\in -\rmd v, S_i\in \rmd x, \sigma_i-\sigma_1\in \rmd s)\\
&={\bf 1}_{\{v>0\}}\sum_{i\ge 2}P(\gb\ge i-1, S_{i-1}\in -\rmd v)F(\rmd x+v)P(\sigma_{i-1}\in \rmd s)\\
&={\bf 1}_{\{v>0\}}\sum_{i\ge 0}P(\gb\ge i, S_{i}\in -\rmd v)P(\sigma_{i}\in \rmd s)F(\rmd x+v)\\
&={\bf 1}_{\{v>0\}}\sum_{i\ge 0}\whU(i, \rmd v)P(\sigma_{i}\in \rmd s)F(\rmd x+v)\\
&={\bf 1}_{\{v>0\}}\whV(\rmd s, \rmd v)F(\rmd x+v)\\
&=\whV(\rmd s, \rmd v)F(\rmd x+v)\\
\ea\een
since $s>0$ and $\whV(\rmd s, \{0\})=\delta_0(\rmd s)$ by \eqref{tV0} (or \eqref{wlp}).\hfill\halmos\bigskip




\bigskip \noindent {\bf Proof of Theorem \ref{CDK},
Compound Poisson case:}\
Since
$X$ does not creep in this case,
we need to prove \eqref{DKcor6} for all $s\ge 0, x\ge 0$.
When $s= x=0$, \eqref{DKcor6} reduces to showing
\be\label{Pi00}
\int_{v\ge 0}\whV(\{0\}, \rmd v)\Pi_X(\{v\})=0.
\ee
But $\whV(\{0\}, \rmd v)= \delta_0(\rmd v)$ by \eqref{tV0},
which proves \eqref{Pi00}.
Thus we may assume $s\vee x>0$.

If $s>0$ and $x\ge 0$, then integrating out $v$ in \eqref{XtV} gives
\be\label{VDe}
\int_{v\ge 0} \whV(\rmd s, \rmd v)F(\rmd x+v)=P(\Delta L^{-1}_{\sigma_1}\in \rmd s,  \Delta H_{\sigma_1}\in \rmd x).
\ee
This continues to hold when  $s=0$ and $x> 0$ since
\ben
\int_{v\ge 0} \whV(\{0\}, \rmd v)F(\rmd x+v)=\int_{v\ge 0} \delta_0(\rmd v)F(\rmd x+v)=F(\rmd x),
\een
while
\ben
P(\Delta L^{-1}_{\sigma_1}=0,  \Delta H_{\sigma_1}\in \rmd x)
=P(X_{\sigma_1}\in \rmd x)
=F(\rmd x).
\een
Thus \eqref{VDe} holds whenever $s\vee x>0$.
Now
by the compensation formula (which requires $s\vee x>0$) we have
\ben\ba
P(\Delta L^{-1}_{\sigma_1}\in \rmd s,  \Delta H_{\sigma_1}\in \rmd x)
&=E\sum_{t\ge 0}{\bf 1}(L^{-1}_{t-}=t, H_{t-}=0,
\Delta L^{-1}_{t}\in \rmd s, \Delta H_t\in \rmd x) \\
&=\int_{t\ge 0} P(L^{-1}_{t}=t, H_{t}=0)\rmd t\ \Pi_{L^{-1},H}(\rmd s, \rmd x) \\
&=\int_{t\ge 0} P(\sigma_1>t)\ \rmd t\ \Pi_{L^{-1},H}(\rmd s, \rmd x)\\
&=\lambda^{-1} \Pi_{L^{-1},H}(\rmd s, \rmd x),
\ea\een
which together with  $\Pi_X(\rmd x)=\lambda F(\rmd x)$ completes the proof of \eqref{DKcor6}.
\hfill\halmos
\bigskip


\end{document}